\documentclass{article}
\usepackage{amssymb}
\usepackage{amsmath}
\usepackage{amsfonts}

\newtheorem{theorem}{Theorem}

\newtheorem{corollary}[theorem]{Corollary}

\newtheorem{definition}[theorem]{Definition}

\newtheorem{remark}[theorem]{Remark}

\newenvironment{proof}[1][Proof]{\noindent\textbf{#1.} }{\ \rule{0.5em}{0.5em}}
\newenvironment{proof1}[1][Proof of Theorem 5]{\noindent\textbf{#1.} }{\ \rule{0.5em}{0.5em}}
\renewcommand{\Re}{\mathop{\rm Re}\nolimits}

\input{tcilatex}

\begin{document}

\title{A generalized version of the Earle-Hamilton fixed point theorem for the Hilbert
ball.}
\author{David Shoikhet}
\maketitle

\begin{abstract}
Let $D$ be a bounded domain in a complex Banach space. According
to the Earle-Hamilton fixed point theorem, if a holomorphic
mapping $F : D \mapsto D$ maps $D$ strictly into itself, then it
has a unique fixed point and its iterates converge to this fixed
point locally uniformly. Now let $\mathcal{B}$ be the open unit
ball in a complex Hilbert space and let $F : \mathcal{B} \mapsto
\mathcal{B}$ be holomorphic. We show that a similar conclusion
holds even if the image $F(\mathcal{B})$ is not strictly inside
$\mathcal{B}$, but is contained in a horosphere internally tangent
to the boundary of $\mathcal{B}$. This geometric condition is
equivalent to the fact that $F$ is asymptotically strongly
nonexpansive with respect to the hyperbolic metric in
$\mathcal{B}$.
\end{abstract}

\maketitle

\section{Introduction}

Let $\Delta $ be the open unit disk in the complex plane $%
\mathbb{C}
$, and let $F$ be a holomorphic self-mapping of $\Delta .$
Combining the classical Denjoy-Wolff Theorem and the
Julia-Wolff-Carath\'{e}odory Theorem (see, \cite{DS}, \cite{C-M},
\cite{G-R}, \cite{SD} and \cite{SJ}), it can be stated that:

\vspace{8pt}

$\blacklozenge $ \textit{If }$F$\textit{\ is not the identity and
is not an elliptic automorphism of }$\Delta ,$\textit{\ then the
iterates }$F^{n}(=F\circ F^{n-1})$\textit{, }$n=1,2,...$\textit{, }$F^{0}=I$\textit{%
, the identity mapping on }$\Delta $\textit{, converge uniformly
on the
compact subsets of }$\Delta $\textit{\ to a constant mapping of }$\Delta $%
\textit{\ into }$\overline{\Delta }$\textit{.}

\textit{Moreover, the following are equivalent:}

\textit{(i) }$F$\textit{\ is fixed point free in }$\Delta
$\textit{\ (that is there is no fixed point of }$F$\textit{\ in
}$\Delta )$.

\textit{(ii) }$F$\textit{\ has a boundary regular fixed point }$\zeta
(=\lim_{r\rightarrow 1^{-}}F(r\zeta ))$\textit{, such that \ }

$\ \ \ \ \ 0<F^{\prime }(\zeta )(=\lim_{r\rightarrow 1^{-}}F^{\prime
}(r\zeta ))\leq 1$\textit{.}

\textit{(iii) There is a point }$\eta \in \partial \Delta $\textit{\ such
that each horodisk internally tangent to }

$\ \ \ \ \ \partial \Delta $\textit{\ at }$%
\eta $ \textit{is }$F-$\textit{invariant.}

\textit{(iv) The sequence of iterates }$\left\{ F^{n}\right\}
_{n=1}^{\infty }$\textit{converges uniformly on the compact }

$\ \ \ \ \ $ \textit{subsets of }$\Delta $\textit{\ to a boundary
point }$\tau \in
\partial \Delta .$

\textit{The points }$\zeta $\textit{, }$\eta $\textit{\ and }$\tau $\textit{%
\ in (ii)-(iv) must coincide.}

\vspace{8pt}

Sometimes this assertion is called the Grand Fixed Point Theorem
(see, for example,\cite{SJ}). In general, it is no longer true for
the infinite dimensional case. Various partial analogs of the
Denjoy-Wolff and Julia-Wolff-Carath\'{e}odory Theorems for higher
dimensions can be found in \cite{AM-98}, \cite{AM-TR-99},
\cite{C-M}, \cite{G-R}, \cite{KT-R-S}, \cite{MBD},
\cite{RS-SD-97b}, \cite{RS-SD-2002a}, \cite{SD}.

Let now $\mathcal{B}$ be the open unit ball in the complex Hilbert space $%
\mathcal{H}$, with the scalar product $\left\langle \cdot ,\cdot
\right\rangle $ and the norm $\left\Vert z\right\Vert =\sqrt{\left\langle
z,z\right\rangle }$, $z\in $ $\mathcal{H}$.

In this paper we study some geometric and analytic conditions which ensure
the local uniform convergence of iterates of a holomorphic self-mapping of $%
\mathcal{B}$ to a constant mapping of the norm less or equal to 1.

We need the following notions and facts.

Fix $a\in\mathcal{B}$. The M\"{o}bius transformation $m_{a}$ on
$\mathcal{B}$ is defined by
\begin{equation}
m_{a}(z)=\frac{1}{1+\left\langle z,a\right\rangle }\left(
\sqrt{1-\left\Vert a\right\Vert^{2} }Q_{a}+P_{a}\right) (z+a)
\label{n1}
\end{equation}%
where $z\in \mathcal{B}$, $P_{a}$ is the orthogonal projection of $\mathcal{H%
}$ onto the subspace $\{\lambda a:\lambda \in
\mathbb{C}
\}$ and $Q_{a}=I-P_{a}$, $I$ is the identity mapping on $\mathcal{H}$.

The Poincar\'{e} hyperbolic metric on $\mathcal{B}$ is the function $\rho _{%
\mathcal{B}}:\mathcal{B\times B\rightarrow
\mathbb{R}
}^{+}$ given by
\begin{eqnarray}
\rho _{\mathcal{B}}(z,w) &=&\tanh ^{-1}\left\Vert m_{-z}(w)\right\Vert =
\label{n2} \\
&=&\frac{1}{2}\log \frac{1+\left\Vert m_{-z}(w)\right\Vert }{1-\left\Vert
m_{-z}(w)\right\Vert },\text{ \ \ }z,w\in \mathcal{B}\text{.}  \notag
\end{eqnarray}

Note that $\rho _{\mathcal{B}}$ can be written in the form%
\begin{equation}
\rho _{\mathcal{B}}(z,w)=\tanh ^{-1}\sqrt{1-\sigma (z,w)}  \label{n3}
\end{equation}%
where%
\begin{eqnarray}
\sigma (z,w) &=&1-\left\Vert m_{-z}(w)\right\Vert ^{2}=  \label{n4} \\
&=&\frac{(1-\left\Vert z\right\Vert ^{2})(1-\left\Vert w\right\Vert ^{2})}{%
\left\vert 1-\left\langle z,w\right\rangle \right\vert ^{2}}\text{.}  \notag
\end{eqnarray}

\begin{definition}
\label{def.1} A mapping $F:\mathcal{B}\rightarrow \mathcal{B}$ is
called $ \rho _{\mathcal{B}}$-nonexpansive if
\begin{equation}
\rho _{\mathcal{B}}(F(z),F(w))\leq \rho _{\mathcal{B}}(z,w)  \label{n5}
\end{equation}%
for all $z,w\in \mathcal{B}$.
\end{definition}

In accordance with (\ref{n3})-(\ref{n5}) $F:\mathcal{B}\rightarrow
\mathcal{B}$\textit{\ is }$ \rho
_{\mathcal{B}}$\textit{-nonexpansive if and only if}
\begin{equation}
\sigma (z,w)\leq \sigma (F(z),F(w))  \label{n6}
\end{equation}%
\textit{for all }$z,w$\textit{\ in }$B$\textit{.}

A consequence of the Schwarz-Pick Lemma is the fact that
\textit{each
holomorphic self-mapping of }$\mathcal{B}$\textit{\ is }$\rho _{\mathcal{B}}$%
\textit{-nonexpansive }(see, for example, \cite{F-V}, \cite{G-R}, \cite{R-S}%
, \cite{DS}).

Since $\mathcal{B}$ is a complete metric space with respect to the metric $%
\rho _{\mathcal{B}}$, the latter fact is very useful in the study
of the fixed point sets of holomorphic self-mappings in
$\mathcal{B}$ in the framework of the general fixed point theory
on metric spaces.

At the same time, since $\rho _{\mathcal{B}}(z,w)$ goes to infinity if
either $z$ or $w$ tends to the boundary $\partial \mathcal{B}$ of $\mathcal{B%
}$ little information can be derived about boundary behavior of
holomorphic (or $\rho _{\mathcal{B}}$-nonexpansive) mappings even
if they admit continuous extension onto $\partial \mathcal{B}$.

To avoid this deficiency, one induces another non-euclidean "distance" $%
d(z,w) $ from $z\in \mathcal{B}$ to $w\in \overline{\mathcal{B}}$, the
closure of $\mathcal{B}$, defined by%
\begin{equation}
d(z,w)=\frac{\left\vert 1-\left\langle z,w\right\rangle \right\vert ^{2}}{%
1-\left\Vert z\right\Vert ^{2}}  \label{n7}
\end{equation}
(see, for example, \cite{G-R}, \cite{C-M}, and \cite{E-R-S2002}).

Note that $d(z,w)$ is not a metric (it is even not symmetric for
$z,w\in \mathcal{B}$, $\| z\| \neq \| w\|$) on $\mathcal{B}$.

Geometrically the sets%
\begin{equation}
E(w,k)=\left\{ z\in \mathcal{B}:d(z,w)<k\right\} ,\text{ }k>1-\left\Vert
w\right\Vert ^{2}\text{, \ }w\in \overline{\mathcal{B}}\text{,}  \notag
\end{equation}%
are ellipsoids in $\overline{\mathcal{B}}$ which for $w\in \partial \mathcal{%
B}$ are usually called \textit{horospheres.}

\begin{definition}
\label{def.2}Let $F$ be a continuous self-mapping of $\mathcal{B}$. A point $%
\zeta \in \partial \mathcal{B}$ is called a fixed point for $F$ if $\zeta
=\lim_{r\rightarrow 1^{-}}F(r\zeta )\in \partial \mathcal{B}$. It is called
a boundary regular fixed point if the radial derivative%
\begin{equation*}
\uparrow F^{\prime }(\zeta ):=\lim_{r\rightarrow 1^{-}}\frac{1-\left\langle
F(r\zeta ),\zeta \right\rangle }{1-r}
\end{equation*}%
exists finitely.
\end{definition}

We use the symbol $\uparrow F^{\prime }(\zeta )$ to distinguish
the radial derivative (which is actually a positive real number,
see \cite{C-M}, \cite{R-S}) from the Frech\'{e}t derivative
usually denoted by $F^{\prime }(z)$ at the point $z\in \mathcal{B}
$ which is a complex linear operator on $\mathcal{H}$ (see, for
example, \cite{F-V}, \cite{G-R}) .

\bigskip

\begin{definition}
\label{def.A} Let $F:\mathcal{B}\rightarrow \mathcal{B}$ be a $\rho _{%
\mathcal{B}}$-nonexpansive mapping on $\mathcal{B}$. A point $\tau \in
\partial \mathcal{B}$ is called a sink point for $F$ if all ellipsoids $%
E(\tau ,k)$, $k>0$ are invariant for $F$.
\end{definition}

It can be shown (cf. \cite{R-S}, Theorems 5.14 and 5.15)
\textit{that a
point \ }$\tau \in \partial B$\textit{\ is a sink point of a }$\rho _{%
\mathcal{B}}$\textit{-nonexpansive mapping }$F$\textit{\ on }$B$\textit{\ if
and only if it is a boundary regular fixed point of }$F$\textit{\ with }$%
\uparrow F^{\prime }(\zeta )$\textit{\ }$\leq 1.$

A result of Goebel, Sekowski and Stachura \cite{G-S-S} (see also
\cite{G-R}, Theorem 25.2) asserts that:

\vspace{8pt}

$\blacklozenge $ \textit{If a }$\rho $\textit{-nonexpansive (holomorphic)
self-mapping }$F:\mathcal{B}\rightarrow \mathcal{B}$\textit{\ is fixed point
free, then there is a unique sink point }$\tau \in \partial \mathcal{B}$ for
$F$.

\vspace{8pt}

In addition,  B. MacCluer showed \cite{MBD}, that \textit{if
}$\dim \mathcal{H}<\infty $\textit{, then iterates of a fixed
point free holomorphic self-mapping }$B$\textit{\ converge
uniformly on compact subsets on }$\mathcal{B}$\textit{\ to a sink
point }$\tau \in \partial \mathcal{B}.$

However, if $\dim \mathcal{H}>1$, then even for holomorphic
mappings, a converse assertion is no longer true: \textit{if
}$F$\textit{\ has a boundary sink point, then }$F$\textit{\ is not
necessarily fixed point free} (see Theorem 25.1 in \cite{G-R} and
examples there in).

Moreover, in contrast with the finite-dimensional case an example
of A. Stachura \cite{SA} shows that for the infinite dimensional
case iterates of a holomorphic self-mapping of $\mathcal{B}$ do
not necessarily converge to a sink point even if it is a unique
boundary regular fixed point of $F$.

On the other hand, one can show (see the proof of Theorem 30.8 in
\cite{G-R}) that if a $\rho _{\mathcal{B}}$-nonexpansive mapping
$F$  in $\mathcal{B}$
is fixed point free and $\tau $ is its sink point then the inequality%
\begin{equation*}
d(F(z),\tau )\leq d\left (\frac{1}{2}(z+F(z)),\tau \right )
\end{equation*}%
provides the \textit{pointwise convergence} (in the norm of
$\mathcal{H}$) of iterates $F^{n}(z)$ to the point $\tau $ for all
$z\in \mathcal{B}$.

Even though the strong and weak convergence of iterates has been
studied very intensively (\cite{G-R2}, \cite{R},see also book
\cite{G-R} and survey \cite{KT-R-S} and references therein),
little is known about local uniform convergence of iterates for
the infinite dimensional case. Some results concerning this issue
for compact or condensing holomorphic mappings can be found in
\cite{K-K-R} and \cite{K-K-R2}.

\begin{definition}
(cf. \cite{F-V} and \cite{R-S}) We say that a subset $K$ of a bounded domain $%
\mathcal{D}$ in a Banach space $X$ is strictly inside $\mathcal{D}$, if it
is bounded away from the boundary of $\mathcal{D}$, i.e., $\inf_{\substack{ %
z\in \mathcal{K}  \\ w\in \partial \mathcal{D}}}\left\Vert z-w\right\Vert
\geq \varepsilon >0$.

A sequence $\{F^{n}\}_{n=1}^{\infty}$ of mappings on $\mathcal{D}$
is said to be locally uniformly convergent on $D$ if it converges
uniformly (in the norm of $X$) on each ball strictly inside
$\mathcal{D}$.
\end{definition}

The famous Earle-Hamilton Theorem \cite{E-H} (see also  \cite{DS},
\cite{G-R}, \cite{HL} and \cite{R-S}) asserts that:

\vspace{8pt}

$\blacklozenge $ \textit{Let }$F$\textit{\ be a holomorphic self-mapping of
a bounded domain }$D$\textit{\ in a complex Banach space }$X. $ \textit{If }$F$%
\textit{\ maps }$D$\textit{\ into a subset }$K$\textit{\ strictly inside }$D$%
\textit{, then }$F$\textit{\ has a unique fixed point }$\zeta \in D$\textit{%
\ and \ iterates }$\{F^{n}\}_{n=1}^{\infty }$\textit{\ converge to }$\zeta $%
\textit{\ locally uniformly in }$D$\textit{.}

\vspace{8pt}

For $%
\mathbb{C}
^{n}$ this result was partially established earlier by M.Herve
\cite{HM1} as a consequence of the Remmert-Stein Theorem.

In the setting of the Hilbert ball we show that a similar
conclusion holds even if the image $F(\mathcal{B})$ is not
strictly inside $\mathcal{B}$, but
belongs to a horosphere in $\overline{\mathcal{B}},$ the closure of $%
\mathcal{B}$.

To formulate our result we need some additional notations.

For a linear operator $A:X\mapsto X$, we denote by $\Sigma (A)$
the spectrum of $A,$ and by $\Sigma _{p}\left( A\right) =\{\lambda
\in \Sigma (A):\lambda $ is an eigen-value of $A\}$ the point
spectrum of $A.$

By $\Sigma _{\partial \Delta }(A)=\{\lambda \in \Sigma (A)$, $\left\vert
\lambda \right\vert =1\}$ we denote the peripheral spectrum of $A,$
whenever, $\Sigma (A)\subseteq \overline{\Delta },$ the closure of the unit
disk.

It is well known that \textit{if a holomorphic self-mapping
}$F$\textit{\ of a bounded domain} $D$ \textit{\ in a Banach
space} $X$\textit{\ has a fixed point }$\zeta \in D$\textit{, then
}$\Sigma (F^{\prime }(\zeta ))\subseteq \overline{\Delta
}$\textit{\ }(see, for example, \cite{VE} and \cite{R-S}).

\begin{theorem}
\label{teoremaN} Let $F$ be a holomorphic mapping on $\mathcal{B}$ which
maps $\mathcal{B}$ into a horosphere $E(\tau ,m)$ for some $\tau \in
\partial \mathcal{B}$ and $0<m<\infty $, i.e.,
\begin{equation}
d(F(z),\tau )< m<\infty \text{, \ }z\in \mathcal{B}\text{.}
\label{n0}
\end{equation}%
The following assertions hold.

(i) $F$ has at most one fixed point in $\mathcal{B}$.

(ii) If $F(\zeta )=\zeta \in \mathcal{B}$ and $A=F^{\prime }(\zeta )$ has
the peripheral spectrum $\Sigma _{\partial \Delta }(A)$ which belongs to the
point spectrum $\Sigma _{p}\left( A\right) $, then the iterates $%
\{F^{n}\}_{n=1}^{\infty }$ converge to the point $\zeta $ in the topology of
local uniform convergence on $\mathcal{B}$.

(iii) If $F$ has a boundary regular fixed point, then it must be
$\tau$.

(iv) $F$ is fixed point free if and only if $\tau (\in \partial
\mathcal{B})$ is a sink point for $F$ and if and only if iterates
$\{F^{n}\}_{n=1}^{\infty }$ converge to the
point $\tau $ in the topology of the local uniform convergence on $\mathcal{B%
}$.
\end{theorem}

\begin{remark}
It follows from the results in \cite{L-Z} (see, also \cite{K-W} and \cite%
{SA-1}) that a power bounded linear operator $A$ satisfies
condition (ii) of the theorem if and only if it is uniformly
precompact (in the norm operator topology). This fact and the
Cauchy inequalities for the Frech\'{e}t derivatives (see, for
example, \cite{F-V}, \cite{G-R} and \cite{R-S}) imply the
following:
\end{remark}

\begin{corollary}
\label{col.A} Let $F$ be a holomorphic mapping on $\mathcal{B}$
which satisfies (\ref{n0}). Then the sequence
$\{F^{n}\}_{n=1}^{\infty}$ of iterates locally uniformly converges
to a constant mapping in $\overline{\mathcal{B}}$ if and only if
it is locally uniformly precompact.

In particular, if $F(\mathcal{B})$ is relatively compact, then the
sequence $\{F^{n}\}_{n=1}^{\infty}$ always converges to a fixed
point of $F$ in $\mathcal{B}$ in the open compact topology on
$\mathcal{B}$.
\end{corollary}

It should be observed that the key point in the proof of the
Earle-Hamilton Theorem is the fact that \textit{each holomorphic self mapping of }$D$%
\textit{\ which maps }$D$\textit{\ into a subset }$K$\textit{\
strictly inside }$D$\textit{\ is a strict contraction with respect
to a hyperbolic metric on }$D.$ So, the Banach Fixed Point
Principle can be applied.

\vspace{8pt}

To prove our theorem we first establish another metric
characterization of holomorphic mappings which satisfy condition
(\ref{n0}).

\begin{definition}
\label{def.N} We say that a mapping $F:\mathcal{B}\rightarrow
\mathcal{B}$ is asymptotically strongly nonexpansive with respect
to a point $\tau \in
\partial \mathcal{B}$ if there is a positive function $p_{F}$ on
$\mathcal{B}\times \mathcal{B}$ such that

(a) \ $0\leq p_{F}(z,w)\leq 1$,\quad $p_{F}(z,w)\not\equiv
0$,\quad $z,w\in \mathcal{B}$;

(b) \ $\lim_{r\rightarrow 1^{-}}\displaystyle\frac{p_{F}(z,r\tau )}{1-r}:=k>0$ \ for all $%
z\in \mathcal{B}$;

and

(c) \ for all $z,w\in \mathcal{B}$,
\begin{equation*}
\sigma (F(z),F(w))\geq \sigma (z,w)(1-p_{F}(z,w))+p_{F}(z,w)(\geq
\sigma (z,w)).
\end{equation*}
\end{definition}

We remark in passing that another class of strongly nonexpansive
mappings in the Hilbert ball was considered in \cite{R1993} and
\cite{K-R}.

\begin{theorem}
\label{teoremaA'} Let $F$ be a holomorphic self-mapping of $\mathcal{B}$,
such that for some $\tau \in \partial \mathcal{B}$, $\lim\limits_{r\rightarrow 1^{-}}F(r\tau )=\tau$ and the radial derivative $%
\uparrow F^{\prime }(\tau )$ exists finitely. Then $F$ is asymptotically
strongly nonexpansive with respect to $\tau $ if and only if there is a
positive number $m<\infty $ such that condition (\ref{n0}) holds.

Moreover,
\begin{equation}
m\geq \frac{2(\uparrow F^{\prime }(\tau ))}{k},  \label{t2'}
\end{equation}%
where $k$ is defined by condition (b) of Definition \ref{def.N}.
\end{theorem}

\begin{proof}
Assume that (\ref{n0}) holds. For $x\in \mathcal{H}$ define
\begin{equation}
S(x)=\Re \left\langle x,\tau \right\rangle +\left\vert \left\langle x,\tau
\right\rangle \right\vert ^{2}-\left\Vert x\right\Vert ^{2}  \label{t3}
\end{equation}
and let $\Pi (=\Pi _{\tau }):=\{x\in \mathcal{H}:S(x)>0\}$ be
Siegel's
domain in $\mathcal{H}$. So, if%
\begin{equation}
C(z)\left(=C_{\tau }(z)\right):=\frac{1}{1-\left\langle z,\tau
\right\rangle }(z+\tau ) \label{t4}
\end{equation}%
is the Caley transformation of $\mathcal{B}$, then $C(\mathcal{B)}=\Pi $ and
$C^{-1}(\Pi )=\mathcal{B}$.

In addition, for each $z\in \mathcal{B}$ we have%
\begin{eqnarray}
S(C(z)) &=&\Re \left\langle C(z),\tau \right\rangle +\left\vert
\left\langle C(z),\tau \right\rangle \right\vert ^{2}-\left\Vert
C(z)\right\Vert ^{2}= \label{t5} \\ &=&\Re \frac{\left\langle
z,\tau \right\rangle +1}{1-\left\langle z,\tau \right\rangle
}+\frac{\left\vert 1+\left\langle z,\tau \right\rangle \right\vert
^{2}}{\left\vert 1-\left\langle z,\tau \right\rangle \right\vert
^{2}}-\frac{\left\Vert z\right\Vert ^{2}+1+2\Re \left\langle
z,\tau \right\rangle }{\left\vert 1-\left\langle z,\tau
\right\rangle \right\vert ^{2}}  \notag \\ &=&\frac{1-\left\Vert
z\right\Vert ^{2}}{\left\vert 1-\left\langle z,\tau
\right\rangle \right\vert ^{2}}=\frac{1}{d( z,\tau ) }%
\text{.}  \notag
\end{eqnarray}%
Therefore, condition (\ref{n0}) can be rewritten in the form%
\begin{equation}
S(C(F(z)))> \frac{1}{m}=:a>0\text{.}  \label{t6}
\end{equation}

On the other hand, since by (\ref{t3})%
\begin{equation}
S(x-a\tau )=\Re \left\langle x,\tau \right\rangle -a+\left\vert \left\langle
x,\tau \right\rangle -a\right\vert ^{2}-\left\Vert x-a\tau \right\Vert
^{2}=S(x)-a\text{,}  \label{t7}
\end{equation}%
it follows from (\ref{t6}) that the mapping $F_{1}$:%
\begin{equation}
F_{1}(z)=C^{-1}(C(F(z))-a\tau )  \label{t8}
\end{equation}%
is well defined holomorphic mapping on $\mathcal{B}$ and maps $\mathcal{B}$
into itself.

Therefore, we get that $F_{1}$ is $\rho _{\mathcal{B}}$-nonexpansive on $%
\mathcal{B}$ which can be expressed by the inequality

\begin{equation}
\sigma (F_{1}(z),F_{1}(w))\geq \sigma (z,w)  \label{t9}
\end{equation}%
for all $z,w\in $ $\mathcal{B}$.

Furthermore, setting $x=C(z)$ and $y=C(w)$ one calculates that%
\begin{equation}
\sigma (z,w)=\frac{4S(x)S(y)}{\left\vert T(x,y)\right\vert ^{2}}\text{,}
\label{t10}
\end{equation}%
where%
\begin{equation}
T(x,y)=\left\langle x,\tau \right\rangle +\left\langle \tau ,y\right\rangle
+2\left( \left\langle x,\tau \right\rangle \left\langle \tau ,y\right\rangle
-\left\langle x,y\right\rangle \right) \text{.}  \label{t11}
\end{equation}%
In addition, it follows from (\ref{t11}) that
\begin{equation}
T(x-a\tau ,y-a\tau )=T(x,y)-2a\text{.}  \label{t12}
\end{equation}

If we set now $x_{1}=x-a\tau $, $x=C(F(z))$ and $y_{1}=y-a\tau $, $y=C(F(w))$
then we have by (\ref{t7})-(\ref{t12})%
\begin{eqnarray}
\sigma (F_{1}(z),F_{1}(w)) &=&\frac{4S(x_{1})\cdot S(y_{1})}{|T(x_{1},y_{1})|^{2}}=\frac{4S(x-a\tau )\cdot S(y-a\tau )}{%
\left\vert T(x-a\tau ,y-a\tau )\right\vert ^{2}} \label{t13} \\
&=&\frac{4(S(x)-a)\cdot (S(y)-a)}{\left\vert T(x,y)-2a\right\vert
^{2}} \notag \\
&=&\frac{4S(x)S(y)-4a[S(x)+S(y)]+4a^{2}}{\left\vert
T(x,y)\right\vert ^{2}-4a\Re T(x,y)+4a^{2}}\text{.}  \notag
\end{eqnarray}

Then by using (\ref{t10}), (\ref{t13}) we obtain%
\begin{equation}
\sigma (F_{1}(z),F_{1}(w))=\frac{\sigma
(F(z),F(w))-b_{F}(z,w)}{1-c_{F}(z,w)}, \label{t15}
\end{equation}%
where we denote
\begin{equation}
b_{F}(z,w)=\frac{4a[S(C(F(z)))+S(C(F(w)))-a]}{\left\vert
T(C(F(z)),C(F(w)))\right\vert ^{2}}  \label{t16}
\end{equation}%
and
\begin{equation}
c_{F}(z,w)=\frac{4a\func{Re}[T(C(F(z)),C(F(w)))-a]}{\left\vert
T(C(F(z)),C(F(w)))\right\vert ^{2}}\text{.}  \label{t17}
\end{equation}

Since $S(C(F(z)))\geq a$, it is clear that $b_{F}(z,w)>0$.

Now we prove that
\begin{equation}\label{nnn1}
0<b_{F}(z,w)\leq c_{F}(z,w)< 1, \quad z,w\in\mathcal{B}.
\end{equation}%

To do this we have to show that%
\begin{equation}
\func{Re}T(x,y)\geq S(x)+S(y)  \label{t18'}
\end{equation}%
and%
\begin{equation}
4a(\func{Re}T(x,y)-a)\leq \left\vert T(x,y)\right\vert
^{2}\text{,} \label{t19'}
\end{equation}%
whenever $\min (S(x),S(y))\geq a$, for $x, y\in \Pi $.

Indeed, by (\ref{t11})%
\begin{equation*}
\Re T(x,y)=\Re [\left\langle x,\tau \right\rangle +\left\langle
\tau ,y\right\rangle +2(\left\langle x,\tau \right\rangle
\left\langle \tau ,y\right\rangle -\left\langle x,y\right\rangle
)]\text{.}
\end{equation*}

Hence, by (\ref{t3}) inequality (\ref{t18'}) can be rewritten as%
\begin{equation*}
2\Re (\left\langle x,\tau \right\rangle \left\langle \tau
,y\right\rangle -\left\langle x,y\right\rangle )\geq \left\vert
\left\langle x,\tau \right\rangle \right\vert ^{2}+\left\vert
\left\langle y,\tau \right\rangle \right\vert ^{2}-\left\Vert
x\right\Vert ^{2}-\left\Vert y\right\Vert ^{2}
\end{equation*}%
or%
\begin{equation}
\left\Vert x-y\right\Vert ^{2}\geq \left\vert \left\langle x-y,\tau
\right\rangle \right\vert ^{2}  \label{t20a}
\end{equation}%
which is obvious.

Note that, in fact, (\ref{t20a}) holds for all $x,y\in
\mathcal{H}$, so does (\ref{t18'}). Inequality (\ref{t19'}) is
also evident since in view of (\ref{t18'}) and (\ref{t6}) we have
that
\begin{equation*}
\left\vert T(x,y)\right\vert ^{2}-4a\Re T(x,y)+4a^{2}=\left\vert
T(x,y)-2a\right\vert ^{2}\geq [S(x)+S(y)-2a]^{2} > 0\text{.}
\end{equation*}
Thus, relation (\ref{nnn1}) is proved.

Furthermore, we calculate by using (\ref{t11}) and (\ref{t4})%
\begin{equation}  \label{n28'}
\begin{split}
T(C(&F(z)),C(F(w))) =\left\langle C(F(z)),\tau \right\rangle +\left\langle
\tau ,C(F(w))\right\rangle {} \\
&+2\left( \left\langle C(F(z)),\tau \right\rangle \left\langle \tau
,C(F(w))\right\rangle - \left\langle C(F(z)),C(F(w))\right\rangle\right ) {}
\\
&=\frac{1}{1-\langle\tau , F(w) \rangle}\cdot\left[\frac{1+\langle
F(z), \tau\rangle}{1-\langle F(z),
\tau\rangle}\cdot\left(1-\langle\tau ,F(w)\rangle \right )
+1+\langle\tau , F(w)\rangle \right . {} \\ & +2
\left(\frac{1+\langle F(z), \tau\rangle}{1-\langle F(z),
\tau\rangle} \left(1+\langle\tau , F(w)\rangle\right) \right. {}
\\ & \left. \left. -\left\langle\frac{1}{1-\langle
F(z),\tau\rangle}(F(z)+\tau) , F(w)+\tau \right\rangle
\right)\right] \text{.} {} \\
\end{split}%
\end{equation}

Consider now the following function%
\begin{equation*}
A(r)=(1-r)^{2}\left\vert T(C(F(z)),C(F(r\tau )))\right\vert
^{2}\text{.}
\end{equation*}

Since
\begin{equation*}
\lim_{r\rightarrow 1^{-}}F(r\tau )=\tau
\end{equation*}%
and
\begin{equation}  \label{t22'}
\lim_{r\rightarrow 1^{-}}\frac{\left\vert 1-\left\langle \tau ,F(r\tau
)\right\rangle \right\vert }{1-r}=\lim_{r\rightarrow 1^{-}}\frac{%
1-\left\Vert F(r\tau )\right\Vert }{1-r}=\uparrow F^{\prime }(\tau
),
\end{equation}
(see Theorems 5.12 - 5.14 in \cite{R-S}) we obtain from
(\ref{n28'})

\begin{equation}
\begin{split}
\lim_{r\rightarrow 1^{-}}A(r)& =\lim_{r\rightarrow
1^{-}}\frac{(1-r)^{2}}{ \left\vert 1-\left\langle \tau ,F(r\tau
)\right\rangle \right\vert ^{2}} \cdot \lim_{r\rightarrow
1^{-}}\left\vert \frac{1+\langle F(z),\tau \rangle }{1-\langle
F(z),\tau \rangle }\cdot \right. {} \\ & \left(1-\langle \tau
,F(r\tau )\rangle \right) +1+\langle \tau ,F(r\tau )\rangle
+2\left( \frac{1+\langle F(z),\tau \rangle }{1-\langle F(z),\tau
\rangle }\cdot \right. {} \\ & \left. \left. \left( 1+\langle \tau
,F(r\tau )\rangle \right) -\left\langle \frac{F(z)+\tau
}{1-\langle F(z),\tau \rangle },F(r\tau )+\tau \right\rangle
\right) \right\vert {} \\ & =4\lim_{r\rightarrow
1^{-}}\frac{(1-r)^{2}}{\left\vert 1-\left\langle \tau ,F(r\tau
)\right\rangle \right\vert ^{2}}=\frac{4}{(\uparrow F^{\prime
}(\tau ))^{2}}\text{.}
\end{split}
\label{t23''}
\end{equation}

Hence,
\begin{equation}
\lim_{r\rightarrow 1^{-}}(1-r)\left\vert T(C(F(z)),C(F(r\tau )))\right\vert
^{2}=\lim_{r\rightarrow 1^{-}}\frac{A(r)}{1-r}=\infty \text{.}  \label{t24}
\end{equation}

Then we obtain from (\ref{t16}) and (\ref{n28'})-(\ref{t24}):
\begin{eqnarray}\label{nnn2}
\lim _{r\rightarrow 1^{-}}\frac{b_{F}(z,r\tau
)}{1-r}&=&\lim_{r\rightarrow 1^{-}}\frac{4a\left[
S(C(F(z)))+S(C(F(r\tau )))-a\right] }{(1-r)\left\vert
T(C(F(z)),C(F(r\tau ))\right\vert ^{2}} \notag
\\ &=&\lim_{r\rightarrow 1^{-}}\frac{4a(1-r)\left[
S(C(F(z)))+S(C(F(r\tau )))-a\right] }{A(r)} \notag \\ &=&a\cdot
\left( \uparrow F^{\prime }(\tau )\right) ^{2}\cdot
\lim_{r\rightarrow 1^{-}}(1-r)S(C(F(r\tau ))) \\ &=&a\cdot
(\uparrow F^{\prime }(\tau ))^{2}\lim_{r\rightarrow
1^{-}}\frac{(1-r)(1-\left\Vert F(r\tau )\right\Vert
^{2})}{\left\vert 1-\left\langle F(r\tau ),\tau \right\rangle
\right\vert ^{2}} \notag \\ &=&2a(\uparrow F^{\prime }(\tau
))=\frac{2(\uparrow F^{\prime }(\tau ))}{m}>0\text{.} \notag
\end{eqnarray}

Since%
\begin{equation*}
\frac{\Re T(x,y)}{\left\vert T(x,y)\right\vert }\leq 1
\end{equation*}%
we also  get from (\ref{t17}) and (\ref{t24})%
\begin{equation}
\begin{split}\label{nnn3}
0 &\leq \lim_{r\rightarrow 1^{-}}c_{F}(z,r\tau )\leq  \\ &\leq
 4a\lim_{r\rightarrow 1^{-}}\left[ \frac{1}{\left\vert
T(C(F(z)),C(F(r\tau )))\right\vert }-\frac{a}{\left\vert
T(C(F(z)),C(F(r\tau )))\right\vert ^{2}}\right] =0.
\end{split}
\end{equation}

To proceed, we need the following observation:

If $\rho _{\mathcal{B}}(\cdot ,\cdot )$ is the Poincar\'{e}
hyperbolic metric on $\mathcal{B}$, then one can define the
hyperbolic metric on $\Pi
=C(\mathcal{B)}$ by the equality%
\begin{equation*}
\rho _{\Pi }(x,y)=\rho _{\mathcal{B}}(C^{-1}(x),C^{-1}(y))\text{,}
\end{equation*}%
whenever, $x$ and $y$ are in $\Pi $. Since the affine mapping $h$ defined by%
\begin{equation*}
h(x)=x+a\tau
\end{equation*}%
maps $\Pi $ into itself we have that%
\begin{equation*}
\rho _{\Pi }(x+a\tau ,y+a\tau )\leq \rho _{\Pi }(x,y)
\end{equation*}%
or%
\begin{equation*}
\rho _{\mathcal{B}}(C^{-1}(x+a\tau ),C^{-1}(y+a\tau ))\leq \rho _{\mathcal{B}%
}(C^{-1}(x),C^{-1}(y))\text{, \ \ \ }x,y\in \Pi \text{.}
\end{equation*}

Since $S(C(F(z))-a\tau )=S(C(F(z)))-a> 0$ and $S(C(F(w))-a\tau
)=S(C(F(w)))-a> 0$ for all $z,w\in \mathcal{B}$, one can set $
x=C(F(z))-a\tau $ and $y=C(F(w))-a\tau $ to obtain
\begin{equation*}
\rho _{\mathcal{B}}(F(z),F(w))\leq \rho _{\mathcal{B}}(F_{1}(z),F_{1}(w))
\end{equation*}%
which is equivalent to%
\begin{equation*}
\sigma (F(z),F(w))\geq \sigma (F_{1}(z),F_{1}(w))=\frac{\sigma
(F(z),F(w))-b_{F}(z,w)}{1-c_{F}(z,w)}\text{.}
\end{equation*}
This implies that
\begin{equation}\label{nnn4}
\sigma (z,w)\leq \sigma (F(z),F(w))\leq \frac{b_{F}(z,w)}{c_{F}(z,w)}\leq 1%
\text{.}
\end{equation}
In turn, the last inequality together with (\ref{t9}) and
(\ref{t15}) implies that
\begin{equation}
\sigma (F(z),F(w))\geq (1-c_{F}(z,w))\sigma (z,w)+b_{F}(z,w)\geq
\sigma (z,w)\text{.}  \label{n0'}
\end{equation}
Now let us define a real function
$p_{F}:\mathcal{B}\times\mathcal{B}\mapsto\mathbb{R}$ by the
formula
\begin{equation}\label{nnn5}
p_{F}(z,w)=\frac{b_{F}(z,w)-\sigma (z,w)c_{F}(z,w)}{1-\sigma
(z,w)}.
\end{equation}
It follows from (\ref{nnn4}) that $$p_{F}(z,w)\geq 0,\quad
z,w\in\mathcal{B},\quad z\neq w .$$ Also, since $$ \lim\limits
_{r\rightarrow 1^{-}}\frac{\sigma (z,r\tau )}{1-r}=\lim\limits
_{r\rightarrow 1^{-}}\frac{\left(1-\|z\| ^{2}
\right)\left(1-r^{2}\right)}{|1-\langle z,r\tau\rangle|^{2}(1-r)}
=\frac{2}{d(z,\tau)}<\infty , $$ we obtain from (\ref{nnn2}) and
(\ref{nnn3})
\begin{equation*}
\begin{split}
k:&=\lim\limits _{r\rightarrow 1^{-}}\frac{p_{F}(z,r\tau )}{1-r}
=\lim\limits _{r\rightarrow
1^{-}}\frac{\frac{1}{1-r}\left(b_{F}(z,r\tau )-\sigma(z,r\tau
)c_{F}(z,r\tau )\right)}{1-\sigma (z,r\tau )} \\ &=\lim\limits
_{r\rightarrow 1^{-}}\frac{b_{F}(z,r\tau )}{1-r}=\frac{2(\uparrow
F' (\tau )) }{m}>0 .
\end{split}
\end{equation*}
Finally, we have from (\ref{nnn5}) the equality
$$\left(1-p_{F}(z,w)\right)\sigma
(z,w)+p_{F}(z,w)=\left(1-c_{F}(z,w)\right)\sigma
(z,w)+b_{F}(z,w),$$ which proves together with (\ref{n0'}) the
inequality $p_{F}(z,w)\leq 1$ and condition (c) of Definition
\ref{def.N}.

This completes the proof of the sufficient part of our theorem.
\end{proof}

To prove the converse part of Theorem \ref{teoremaA'} we establish
a more general assertion which holds actually for asymptotically
strongly nonexpansive mappings which are not necessarily
holomorphic. This assertion also proves assertions (iii) and (iv)
of Theorem \ref {teoremaN}.

\begin{theorem}
\label{teor.10''} Let $F:\mathcal{B}\rightarrow \mathcal{B}$ be an
asymptotically strongly nonexpansive with respect to a point $\tau
\in \partial \mathcal{B}$ (see Definition \ref{def.N}). Assume
that $\tau$ is a boundary regular fixed point of $F$, i.e.,
$\lim\limits_{r\rightarrow 1^{-}}F(r\tau )=\tau$ and the radial
derivative $\uparrow F^{\prime }(\tau )=\beta
>0$ exists finitely. Then for all $z\in
\mathcal{B}$ the following inequality holds%
\begin{equation*}
\frac{1}{d(F(z),\tau )}\geq \frac{1}{\beta }\left( \frac{k}{2}+\frac{1}{%
d(z,\tau )}\right)
\end{equation*}%
where $\displaystyle k=\lim_{r\rightarrow
1^{-}}\frac{p_{F}(z,r\tau )}{1-r}$.

If $F$ has a boundary regular fixed point then it must be $\tau$.
In particular, the following are equivalent:

(i) The mapping $F$ is fixed point free.

(ii) The point $\tau$ is a sink point for $F$.

(iii) The sequence of iterates $\left\{ F^{n}\right\}
_{n=1}^{\infty }$ converges to the point $\tau $ locally uniformly
on $\mathcal{B}$.

Moreover, the following rate of convergence holds
\begin{equation*}
d(F^{n}(z),\tau )\leq \alpha (n,z)\cdot d(z,\tau )\text{,}
\end{equation*}%
where%
\begin{equation*}
\alpha (n,z)=\left\{
\begin{array}{c}
\displaystyle\frac{2}{2+nkd(z,\tau )}\text{, if }\beta =1 \\
\\
\displaystyle\frac{2\beta ^{n}}{2(1-\beta )+(1-\beta ^{n})k}\text{, if }\beta <1%
\end{array}%
\right . \text{.}
\end{equation*}
\end{theorem}

\begin{proof}
Let $F$ be a holomorphic self-mapping which satisfies conditions
(a)-(c) of Definition \ref{def.N}.

Then we have%
\begin{eqnarray*}
\frac{(1-\left\Vert F(z)\right\Vert ^{2})(1-\left\Vert F(r\tau )\right\Vert
^{2})}{\left\vert 1-\left\langle F(z),F(r\tau )\right\rangle \right\vert ^{2}%
} &\geq & \\ &\geq &(1-p_{F}(z,r\tau ))\frac{\left( 1-\left\Vert
z\right\Vert ^{2}\right)
(1-r^{2})}{\left\vert 1-\left\langle z,r\tau \right\rangle \right\vert ^{2}}%
+p_{F}(z,r\tau )
\end{eqnarray*}%
or%
\begin{eqnarray*}
\frac{1-\left\Vert F(z)\right\Vert ^{2}}{\left\vert 1-\left\langle
F(z),F(r\tau )\right\rangle \right\vert ^{2}}\cdot
\frac{1-\left\Vert F(r\tau )\right\Vert ^{2}}{1-r^{2}} &\geq & \\
&\geq &(1-p_{F}(z,r\tau ))\frac{1-\left\Vert z\right\Vert
^{2}}{\left\vert
1-\left\langle z,r\tau \right\rangle \right\vert ^{2}}+\frac{p_{F}(z,r\tau )%
}{1-r^{2}}\text{.}
\end{eqnarray*}

Letting $r\rightarrow 1^{-}$ we get%
\begin{equation}\label{mmmm}
\frac{1}{d(F(z),\tau )}\geq \frac{1}{\beta }\left( \frac{k}{2}+\frac{1}{%
d(z,\tau )}\right) .
\end{equation}

Let us assume, that there is a boundary regular fixed point $\eta
$ of $F$, i.e., the radial derivative $L=\uparrow F^{\prime }(\eta
)$ exists finitely.

Then we have by Theorem 5.12 in \cite{R-S}%
\begin{equation*}
d(F(z),\eta )=\frac{\left\vert 1-\left\langle F(z),\eta
\right\rangle \right\vert ^{2}}{1-\left\Vert F(z)\right\Vert
^{2}}\leq L\frac{\left\vert 1-\left\langle z,\eta \right\rangle
\right\vert ^{2}}{1-\left\Vert z\right\Vert ^{2}}=Ld(z,\eta
)\text{.}
\end{equation*}

Choose any $\varepsilon >0$ and element $z\in \mathcal{B}$ such that $%
Ld(z,\eta )<\varepsilon ^{2}$. Setting $w=F(z)$ we have that
\begin{equation*}
\left\vert 1-\left\langle w,\eta \right\rangle \right\vert ^{2}<
\varepsilon ^{2}
\end{equation*}%
and%
\begin{equation}
1-\left\Vert w\right\Vert ^{2}\leq 1-\left\vert \left\langle
w,\eta \right\rangle \right\vert ^{2}\leq 2\left\vert
1-\left\langle w,\eta \right\rangle \right\vert < 2\varepsilon
\text{.}  \label{Nd}
\end{equation}

Also, note that%
\begin{equation}
1-\Re \left\langle w,\eta \right\rangle <\left\vert 1-\left\langle
w,\eta \right\rangle \right\vert < \varepsilon \text{.} \label{Ne}
\end{equation}

On the other hand, since $d(w,\tau )<m$, inequality (\ref{Nd})
implies that
\begin{equation*}
\left\vert 1-\left\langle w,\tau \right\rangle \right\vert ^{2}<
m(1-\left\Vert w\right\Vert ^{2})\leq 2m\varepsilon \text{.}
\end{equation*}

Hence, again%
\begin{equation*}
1-\Re \left\langle w,\tau \right\rangle \leq \left\vert
1-\left\langle
w,\tau \right\rangle \right\vert \leq \delta =\sqrt{2m\varepsilon }%
\text{.}
\end{equation*}

Then we obtain%
\begin{eqnarray*}
\left\Vert \eta -\tau \right\Vert ^{2} &\leq &\left( \left\Vert
\eta -w\right\Vert +\left\Vert \tau -w\right\Vert \right) ^{2}\leq
\\ &\leq &2\left( \left\Vert \eta -w\right\Vert ^{2}+\left\Vert
\tau -w\right\Vert ^{2}\right) = \\ &=&4\left( 1-\Re \left(
\left\langle \eta ,w\right\rangle +\left\langle \tau
,w\right\rangle \right) +\left\Vert w\right\Vert ^{2}\right)  \\
&<&4\left( 1-\Re \left\langle \eta ,w\right\rangle +1-\Re
\left\langle \tau ,w\right\rangle \right)  \\ &< &4\left(
\varepsilon +\delta \right) \text{.}
\end{eqnarray*}

Since $\delta \rightarrow 0^{+}$ as $\varepsilon \rightarrow
0^{+}$ we have that $\tau =\eta $.

It is now clear that if $F$ is fixed point free, then a sink point
$\eta \in\partial\mathcal{B}$ for $F$ (which is a boundary regular
fixed point) must be $\tau$.

Let now $\tau \in \partial \mathcal{B}$ be a sink point for $F.$
In this case $\beta \leq 1$ and we have by induction from
inequality (\ref{mmmm}) that
\begin{eqnarray*}
\frac{1}{d(F^{n}(z),\tau )} &\geq &\frac{1}{\beta ^{n}}\frac{1}{d(z,\tau )}+%
\frac{k}{2\beta }\left( 1+\frac{1}{\beta }+...+\frac{1}{\beta
^{n-1}}\right) =
\\
&=&\left\{
\begin{array}{c}
\displaystyle\frac{1}{d(z,\tau )}+\displaystyle\frac{nk}{2}\text{,
\ if }\beta =F^{\prime }(\tau )=1
\\
\\
\displaystyle\frac{1}{\beta ^{n}}\cdot \left[ \displaystyle\frac{1}{d(z,\tau )}+\displaystyle\frac{1-\beta ^{n}}{%
1-\beta }\displaystyle\frac{k}{2}\right] \text{, \ \ \ if }\beta =F^{\prime }(\tau )<1%
\text{.}%
\end{array}%
\right .
\end{eqnarray*}%
This proves the implication (ii)$\Rightarrow$(iii). The
implication (iii)$\Rightarrow$(i) is obvious, and we are done.
\end{proof}

\begin{proof1}
(i) If $F$ has more than one fixed point in $%
\mathcal{B}$, then it follows from Rudin's Theorem (see \cite{RW}), that the
fixed point set $\mathcal{F}=Fix_{\mathcal{B}}F$ is an affine subset of $%
\mathcal{B}$, hence there is $x\in \mathcal{F}$, such that $x(=F(x))\notin
E(\tau ,m)$. A contradiction.

(ii) Suppose that $\mathcal{F}\neq \varnothing $, i.e., there is $\zeta \in
\mathcal{B}$ such that $\zeta =F(\zeta )$.

Assume also that $A=F^{\prime }(\zeta )$ satisfies the condition%
\begin{equation}
\Sigma _{\partial \Delta }(A)\subseteq \Sigma _{p}(A)\text{.}  \label{Na}
\end{equation}%
First we show that, in fact, under our assumptions the peripheral spectrum $%
\Sigma _{\partial \Delta }(A)$ is empty.

Indeed, assume on the contrary that for some $\ominus \in \lbrack
0,2\pi ]$
there is $x\neq 0$, such that%
\begin{equation}
e^{i\ominus }x=Ax . \label{Nb}
\end{equation}

Consider a holomorphic self-mapping $F_{1}$ of $\mathcal{B}$ defined as
follows%
\begin{equation}
F_{1}=e^{-i\ominus }\Phi \circ F\circ \Phi ^{-1}  \label{Nc} ,
\end{equation}%
where $\Phi :=m_{-\zeta}$ is the M\"{o}bius transformation of
$\mathcal{B}$ defined by (\ref{n1}), taking the point $\zeta $ to
the origin $(\Phi (\zeta )=0)$. Then $F_{1}(0)=0$ and, by
the chain rule%
\begin{equation*}
A_{1}:=F_{1}^{\prime }(0)=e^{-i\ominus }B\circ A\circ B^{-1} ,
\end{equation*}%
where $B=\Phi ^{\prime }(\zeta )$.

Now if $x\neq 0$, $x\in \mathcal{B}$ satisfies (\ref{Nb}), then $A_{1}z=z$,
where $z=Bx\in \mathcal{B}$ and $z\neq 0$.

Now again it follows from Rudin's Theorem that
\begin{equation*}
\mathcal{F}_{1}:=Fix_{\mathcal{B}}(F_{1})=\ker (I-A_{1})\cap \mathcal{B}=%
\mathcal{N\neq }\{0\}\text{.}
\end{equation*}

Therefore, one can find $z\in \mathcal{N}$, such that $z\notin \Phi (E(\tau
,m))$, hence $u=\Phi ^{-1}(z)\notin E(\tau ,m)$.

We claim that there is $x\in \mathcal{B}$ such that $F(x)=u$. Indeed, set $%
y=e^{-i\ominus }z$ and $x=\Phi ^{-1}(y)\in \mathcal{B}$. Since $z\in
\mathcal{N}$ the element $y$ also belongs to $\mathcal{N}$ and we have $%
F_{1}(y)=y$. Then we get by (\ref{Nc}) that
\begin{equation*}
u=\Phi ^{-1}(z)=\Phi ^{-1}(e^{i\ominus }y)=\Phi ^{-1}(e^{i\ominus
}F_{1}(y))=F(\Phi ^{-1}(y))=F(x)\text{.}
\end{equation*}%
This contradicts our assumption.

So, we have proved that the whole spectrum $\Sigma (A)$ belongs to $\Delta $.

Since $\Sigma (A)$ is compact, one can find $\varepsilon >0$, such that the
spectral radius $r(A)<1-\varepsilon $.

This means that there is a norm $\left\Vert \cdot \right\Vert _{1}$ in $%
\mathcal{H}$ equivalent to the original norm of $\mathcal{H}$ such that $%
\sup_{\left\Vert x\right\Vert _{1}\leq 1}\left\Vert Ax\right\Vert _{1}<1$
(see, for example, \cite{Kr}). Since $A=F^{\prime }(\zeta )$ and $\zeta
=F(\zeta )$, this implies by the local uniform continuity of $F^{\prime }(x)$
that there are positive numbers $r>0$ and $\delta >0$ small enough such that%
\begin{equation*}
\left\Vert F(x)-\zeta \right\Vert \leq (1-\delta )\left\Vert x-\zeta
\right\Vert
\end{equation*}%
whenever $\left\Vert x-\zeta \right\Vert \leq r$.

Therefore, $\{F^{n}\}_{n=1}^{\infty }$ converge uniformly on the ball $
\{ x: \left \Vert x-\zeta \right\Vert <r\}\subset \mathcal{B}$ to the point $%
\zeta $. Applying Vitali's convergence Theorem (see, for example \cite{R-S})
we complete assertion (ii).

Assertions (iii) and (iv) are direct results of Theorems
\ref{teoremaA'} and \ref{teor.10''}.
\end{proof1}

\vspace{8pt}

\textbf{Acknowledgements}. The author is very thankful to Prof.
Simeon Reich, Prof. Mark Elin, and Dr. Felix Kerdman for useful
remarks and discussions.

\bigskip

\bigskip

\bigskip

\bigskip

\end{document}